\documentclass[reqno,10pt]{article}
\usepackage{amsmath}
\usepackage{bm}
\usepackage{a4wide,amssymb}
\usepackage{graphicx,xcolor}
\usepackage{epstopdf}
\usepackage{bbm}
\usepackage{mathrsfs}
\usepackage[textheight=185mm,textwidth=125mm]{geometry}

\usepackage[titletoc,toc,title]{appendix}

\usepackage[refpage]{nomencl}
\usepackage{etoolbox}

\newcommand{\qed}{ $\blacksquare$} \newcommand{\nn}{\nonumber}
\newcommand{\R}{{\mathbb R}}         
 
\newcommand{\E}{{\mathbb E}}

\newcommand{\g}{\gamma}\renewcommand{\d}{\delta}

\newcommand{\e}{\varepsilon}

\newcommand{\pa}{\partial}

\newcommand{\D}{\Delta}

\def\G{\Gamma}
\def\s{\sigma}
\def\r{{\rho}}

\def\L{{\Lambda}}
\def\II{{\mathcal I}}

\def\NN{{\cal N}}
\def\II{{\cal I}}

\def\RRR{{\cal R}}

\def\AA{{\cal A}}

\def\SS{{\cal S}}
\def\WW{{\cal W}}

\newtheorem{thm}{Theorem}[section]
\newtheorem{prop}[thm]{Proposition}

\def\be{\begin{equation}}
\def\ee{\end{equation}}
\def\bea{\begin{eqnarray}}
\def\eea{\end{eqnarray}}

\def\nn{\nonumber}

\def\d{\delta}

\def\b{\beta}
\def\t{\tau}

\usepackage{multicol}
\makenomenclature
\usepackage{xpatch}
\makeatletter
\xapptocmd\thenomenclature{\let\@item\nomencl@item\def\nomencl@width{0pt}}{}{}
\let\nomencl@item\@item
\xpretocmd\nomencl@item{\nomencl@measure{#1}}{}{}
\def\nomencl@measure#1{%
  \sbox0{#1}%
  \ifdim\wd0>\nomencl@width\relax
    \edef\nomencl@width{\the\wd0}%
  \fi
}

\xpatchcmd\thenomenclature
  {\section*{\nomname}}
  {\begin{multicols}{2}[\section*{\nomname}]}
  {}{}
\xpatchcmd\thenomenclature
  {\chapter*{\nomname}}
  {\begin{multicols}{2}[\chapter*{\nomname}]}
  {}{}

\xapptocmd\endthenomenclature{%
  \immediate\write\@mainaux{\global\nomlabelwidth\nomencl@width\relax}%
  \end{multicols}
}{}{}
\makeatother

\begin{document}

\begin{titlepage}

\begin{center} 
{\bf \Large{Kinetic SIR equations and particle limits}}

\vspace{1.5cm}
{\large A. Ciallella$^{(*)1}$, M. Pulvirenti$^{2}$ and S. Simonella$^{3}$  }
\vspace{0.5cm}

{$1.$ \scshape {\small Dipartimento di Ingegneria Civile, Edile -- Architettura e
             Ambientale, and 
             \\International Research Center M\&MOCS, 
             Universit\'a dell'Aquila,\\
             via Giovanni Gronchi 18, 67100, L'Aquila, Italy,
            \\ email {alessandro.ciallella@univaq.it} .}
\smallskip

$2.$ {\small Dipartimento di Matematica, Universit\`a di Roma La Sapienza\\ 
Piazzale Aldo Moro 5, 00185 Rome -- Italy, \\
 International Research Center M\&MOCS, Universit\`a dell'Aquila, Italy, and\\
 Accademia Nazionale dei Lincei,
            \\ email {pulviren@mat.uniroma1.it} .}  
 \smallskip

$3.$ {\small UMPA UMR 5669 CNRS, ENS de Lyon \\ 46 all\'{e}e d'Italie,
69364 Lyon Cedex 07 -- France,
            \\ email {sergio.simonella@ens-lyon.fr} .}}
\bigskip

{$(*)$ \scshape {\small corresponding author}}

\end{center}

\vspace{1.5cm}
\noindent
{\bf Abstract.} We present and analyze two simple $N$-particle particle systems for the spread of an infection, respectively with binary and with multi-body interactions.  We establish a convergence result, as $N \to \infty$, to a set of kinetic equations, providing a mathematical justification of related numerical schemes.
We analyze rigorously the time asymptotics of these equations, and compare the models numerically.

\thispagestyle{empty}

\bigskip
\bigskip

\newpage

\end{titlepage}

\section{Introduction} \label{sec:1}
\setcounter{equation}{0}    
\def\theequation{1.\arabic{equation}}

The mathematical models for epidemic spread describe the evolution of average fractions for several species of interacting agents, as susceptible-infected-recovered in the classical SIR model \cite{KMK}; see for instance \cite{AM79, Mu02, HS, BC-C01}.  It is hard to go far beyond this description incorporating spatial patterns in a realistic way, although this can play an important role in applications.
A natural approach is based on modelling equations inspired from the kinetic theory of rarefied gases, see e.g.\,\cite {Betal} and references therein. Of course, the main obstacle is the identification of the assumptions characterizing the interaction. Think for instance of individual strategies, which, if important, could lead to a rather different behaviour from that of a particle system with local equilibrium. 
On the other hand, the essential features of the evolution have little dependence on microscopic details: even a naive model based on a three-species Boltzmann equation can capture equally well the qualitative behaviour of SIR-like equations \cite{PS_SIR,CPS20}.

The present paper is devoted to a mathematical analysis of what could be the simplest possible kinetic model. Agents move independently according to a random flight.  
An infected and a susceptible particle can react into a pair of infected particles, whenever they are sufficiently close. Furthermore each infected particle decays into a recovered particle in a random time of order one  (and cannot be infected anymore, as in the SIR model). 

We then consider two possibilities corresponding to two different scalings. In the first one, the infection reaction of a pair happens in a random time of order $N$, where $N$ is the total number of agents.  
In the second one,  in a random time of order $1$,  the infection is instantaneously transmitted to {\em all} susceptible agents which happen to be close enough (``crowd contagion'', or ``superspread event''). 
As we shall see, both models are conceived to give, {\it formally}, the same kinetic equations  in the limit $N \to \infty$.
But the second one is mathematically more involved. 
 
The plan of the paper is the following. In the next section we present the models and the limiting kinetic equation. In Section \ref{sec:3} we study the time asymptotics of these equations and compare it with the behaviour of the standard SIR. In Section \ref{sec:4} we discuss the convergence of the first particle model in the kinetic limit. Finally in Section \ref{sec:5}, we present some numerical simulation.

\section{Models} \label{sec:2}
\setcounter{equation}{0}    
\def\theequation{2.\arabic{equation}}

\subsection{Model 1}

\subsubsection{Phase space, generators}

Consider a set of $N$ particles (agents) whose position and velocity are in the two-dimensional torus $\L=(0,D)^2$, $D>0$ 
 and the unit circle ${\mathbb S}$ respectively. We denote by $|\L|=D^2$ the measure of $\L$. Particle $i\in \{1,\cdots,N\}$ has a label $a_i \in \{S,I,R\}=:L $.  We set $Z_N=(z_i)_{i=1}^N$ with $z_i=(x_i,v_i) \in \L \times {\mathbb S} $ and $A_N=\{a_i \}_{i=1}^N$. A state of the system $(Z_N; A_N)$ lives in the phase space $(\G \times L)^N$, where $\G=\L \times {\mathbb S}$.

The time evolution is given by the generator
\be
\label{gen1}
 {\cal L}= {\cal L}_0+{\cal L}_1+ {\cal L}_d +{\cal L}^N_{int}
\ee
where: 
\begin{itemize}
\item $ {\cal L}_0=\sum_{i=1}^N v_i \cdot \nabla_{x_i}$ describes free  motion.
\item ${\cal L}_1$ describes velocity jumps
$$
{\cal L}_1\Phi (Z_N) =\sum_{i=1}^N \frac 1{2\pi} \int_{{\mathbb S}} dw \, [\Phi (z_1, \dots, x_i,w, \dots, z_N)-\Phi (Z_N)]
$$
and, for a function $\Phi(Z_N) \equiv\varphi (z_1)$,
$$
{\cal L}_1\Phi (Z_N) =\frac 1{2\pi} \int_{{\mathbb S}} dw  \,[\varphi  (x_1,w)-\varphi (z_1)]\;.
$$
Hence ${\cal L}_0+{\cal L}_1$  generates $N$ independent copies of a random flight.
We do not make explicit here the dependence on the labels, which are not involved.
\item ${\cal L}_d$ describes the decay of infected ($I$) into recovered ($R$) particles
\begin{equation*}
 {\cal L}_d\Phi ( Z_N;A_N)= \g \sum_{i=1} ^N 
[ \, \Phi ( Z_N; a_1, \ldots, \tilde a_i, \ldots, a_N) \! - \Phi ( Z_N; a_1, \ldots, a_i, \ldots, a_N )]
\end{equation*}
where $\g>0$ and the transition $a_i \to \tilde a_i$ is defined by
$$
\tilde a_i=R \,\,\,\,\,\, \text {if} \,\,\,\,\,a_i=I; \qquad \tilde a_i=a_i  \,\,\,\,\, \text{otherwise}.
$$
\item
The interaction or ``infection'', acting over the labels of the agents, is described by 
\begin{equation*}
\begin{split}
 {\cal L}^N_{int}\Phi  (Z_N;A_N )= \frac{\lambda}{N} 
 \sum_{\substack {i \leq N, j \leq N \\ i<j}}
  \big[\,\,& \Phi( Z_N; a_1, \dots, a_i', \dots, a_j', \dots, a_N)\\ - &\Phi( Z_N; a_1, \ldots, a_N) \big]
\end{split}
\end{equation*}
where $\lambda>0$ and the transition $(a_i,a_j) \to (a_i', a_j')$ is defined by
\begin{equation}
\label{a'}
\begin{cases}a'_i=a'_j=I \,\,\,
\text{if} \,\,\, a_i=I, a_j =S \,\,\, \text {or} \,\,\, a_j=I, a_i =S  \,\,\, \text{and} \,\,\, \chi_{i,j}=1, \,\,\,  \\
 a_i'=a_i, \, a_j =a'_j \,\,\, \text {otherwise} \\
\end{cases}\, ,\nn
\end{equation}
with $\chi_{i,j}$ the characteristic function of two particles being at distance less than $R_0>0$
$$
\chi_{i,j} = \mathbbm{1}_{\{ x_i,x_j\  |\  |x_i-x_j | <R_0 \}}\;.
$$
\end{itemize}
Note that, to simplify notation, we dropped the $N$ dependence on generators acting as sums over single particle variables (describing independent particles).

In words, we have the following behaviour. $N$ agents of type $S$ (susceptible), $I$ (infected) or $R$ (recovered) are moving in $\L$ via a random flight, with velocity jumps in ${\mathbb S}$ taking place with rate $1$ per agent. Each infected agent becomes recovered according to a Poisson process of rate $\g$.  The infection is also a Poisson process: we choose a pair of agents $ i,j$ with uniform probability and, if the pair is constituted by an infected and a susceptible and if their distance is smaller than $R_0$, both of agents become instantaneously infected (otherwise nothing happens). The intensity of this process scales like $\lambda N$.

\subsubsection{Densities and marginal distributions}

A statistical description is provided in terms of an initial probability density
$$
W_0^N : (\G \times L)^N \to {\mathbb R}^+ \;,
$$
symmetric in the exchange of particles and normalized by
$$
\sum_{A_N} \int dZ_N W_0^N (Z_N;A_N)=1\;.
$$
The time evolved measure $ W_t^N (Z_N;A_N), t>0$ is given by
\begin{equation*}
\begin{split}
&\sum_{A_N} \int dZ_N W_t^N (Z_N;A_N) \Phi (Z_N;A_N)= \\
&\sum_{A_N} \int dZ_N W_0^N (Z_N;A_N)  \E  [\Phi (Z_N(t);A_N(t) )]\;,
\end{split}
\end{equation*}
where $\Phi$ is a test function, $(Z_N;A_N)  \to (Z_N(t);A_N(t))$ is the process and $\E=\E_{(Z_N,A_N)} $ is the expectation conditioned to the initial value $(Z_N;A_N)$.
The $j$-particle marginal, $j=1, \ldots, N$,  is defined by
$$
f^N_j(Z_j;A_j;t)= \sum_{A'_{N-j} \in L^{N-j}} \int dZ'_{N-j} W_t^N ( Z_j,Z'_{N-j};A_j,A'_{N-j})\;;
$$
giving the probability density of finding $j$ agents with labels $A_j$ in the configuration $Z_j$.

We shall assume full independence at time zero:
\begin{equation} \label{eq:fullindtz}
W_0^N (Z_N;A_N)=\prod_{i=1}^N f_0(x_i,v_i; a_i)
\end{equation}
 where $f_0$ is a one-particle density distribution, with normalization
\begin{equation} \label{eq:norm0}
\sum_{a\in L} \int dz f_0(z;a)=1\;.
\end{equation}
As usual in kinetic theory, the dynamics creates correlations and the measure is not factorized anymore at positive times. One hopes to recover such independence in the limit $N \to \infty$, thanks to the mean-field nature of the interaction (``propagation of chaos'').
Indeed the agents move independently and, given a pair of particles, say $1$ and $2$, the probability that the label of $2$ influences the label of $1$ is $O(\frac 1N)$; therefore that any $j$-particle marginals factorize, if they do factorize at time zero as guaranteed by \eqref{eq:fullindtz} (``propagation of chaos'').

\subsubsection{Kinetic limit} \label{sec:KL}

Let us derive formally the kinetic equations in the limit $N\to \infty$. We choose a test function of type $\Phi(Z_N;A_N)= \frac 1N \sum_{i=1}^N \phi (z_i;a_i) $. 
Then by using the symmetry
\bea
\frac d{dt} \sum_a \int dz  f_1^N \phi (z;a) &=& \frac d{dt} \sum_{A_N} \int dZ_N W^N \Phi (Z_N;A_N)\nn \\
&=& \sum_a \int dz f_1^N( {\cal L}_0 +{\cal L}_1+ {\cal L}_d ) \phi (z;a) \nn \\
&& + \sum_{A_N} \int dZ_N W^N   {\cal L}^N_{int} \Phi (Z_N;A_N)\;, \nn
\eea
where the last term reads
\begin{equation*}
\begin{split}
\label{ha}
  &\frac \lambda {N^2} \sum_{\substack {i,j \\ i < j} }   \sum_{A_N}  \int dZ_N W^N(Z_N;A_N) 
[ \phi (z_i;a'_i) - \phi (z_i;a_i) +\phi (z_j;a'_j) - \phi (z_j;a_j) ] = \nn \\
 & \lambda  \frac {N-1} {2N}\!\! \sum_{ a_1,a_2}\!\! \int \!\!dz_1 dz_2 f^N_2 (z_1,z_2;a_1, a_2) 
[ \phi (z_1;a'_1)\! - \phi (z_1; a_1) 
 +\phi (z_2;a'_2)\! - \phi (z_2;a_2) ] \, . 
\end{split}
\end{equation*}
If $\phi(z;a)=0$ for $a \neq R$, then the interaction term is vanishing. 
If instead $\phi(z;a)=0$ for $a \neq S$,  
the interaction term is close to
$$
 -\lambda \int dz_1  dz_2 f(z_1; S) f(z_2;I)  \,\chi_{1,2}\, \phi (z_1;S)
$$
for $N$ large, if the propagation of chaos holds: $f^N_2 \approx (f_1^N)^{\otimes 2}$ and $\displaystyle f := \lim_{N\to \infty} f^N_1$.
Similarly if $\phi(z;a)=0$ for $a \neq I$ one gets 
$$
+\lambda \int dz_1 dz_2 f(z_1; I) f(z_2;S) \, \chi_{1,2}\, \phi (z_2;I)\;.
$$

We conclude that the triple $ \left(f(z;S;t), f(z;I;t), f(z;R;t)\right)$ satisfies the following system of  kinetic equations ($z=(x,v)$):
\begin{equation}
\label{B}
\begin{cases}
  \left(\partial_t + v \cdot \nabla_x\right) f(z;S) =&  {\cal L}_1 f (z;S) - \lambda  f(z;S) \int f(z_1;I) \chi (|x-x_1| <R_0) dz_1 \\
    \left(\partial_t + v \cdot \nabla_x\right) f(z;I)  =& {\cal L}_1f(z;I)  -\g f(z;I)
    \\ &\qquad+\,\,\lambda  f(z;S)  \int f(z_1 ;I) \chi (|x-x_1| <R_0) dz_1   \\
     \left(\partial_t + v \cdot \nabla_x\right) f (z;R) =& {\cal L}_1 f(z;R ) + \g f(z;I) 
\end{cases}\, .
\end{equation}
Note that the sum
$$
f(z,t) := \sum_{a \in L} f(z;a;t)
$$
satisfies the simple random flight equation
\be
\label{rf}
\left(\partial_t + v \cdot \nabla_x\right) f(z,t) =  {\cal L}_1 f (z,t)\;.
\ee

 \subsection{Model 2}
 \label{ssect:mod2}
 With same setting and notations as above, we now consider the stochastic process with generator
 $$
 \tilde {\cal L}= {\cal L}_0+{\cal L}_1+ {\cal L}_d +\tilde {\cal L}^{N}_{int}
$$
where 
$$
 \tilde {\cal L}^{N}_{int} \Phi (Z_N;A_N)=\frac \lambda N \sum_{i=1}^N 
   \big[\Phi(Z_N;A'_N)- \Phi(Z_N;A_N)\big] \;,
$$
and
$
A'_N=A'_N(i)= \left(a'_j(i) \right)_{j=1}^N
$ is given by
\begin{equation*}
\label{A'}
\begin{cases}
   a'_j =I \,\,\, \text {if} \,\,\, a_j=S, a_i=I  \,\,\, \text {and} \,\,\, \chi_{i,j}=1\\
 a'_j =a_j \,\,\, \text {otherwise} \\
\end{cases}\;.
\end{equation*}
As before, $N$ agents evolve via a random flight and each infected agent becomes recovered according to a Poisson process of rate $\g$; but the spread of the infection affects, with rate $\lambda$, all susceptible agents inside a ball of radius $R_0$ around the infected one.

Proceeding as in Subsection \ref{sec:KL}, we obtain a formal limit by computing $\sum \int W^N \tilde {\cal L}^{N}_{int} \Phi$ for a test function  $\Phi(Z_N;A_N)= \frac 1N \sum_{s=1}^N \phi (z_s;a_s) $:
\begin{equation*}
\begin{split}
&\sum \int W^N \tilde {\cal L}^{N}_{int} \Phi=\\&= \frac \lambda {N^2} \sum_{i} \sum_{s\neq i} \sum_{A_N} \int dZ_N W^N(Z_N;A_N) \chi_{i,s}[ \phi (z_s;a'_s(i)) - \phi (z_s;a_s(i)) ]  \\
&= \lambda \frac {N-1}N \sum_{a_1,a_2} \int dz_1 dz_2 f_2^N (z_1, z_2; a_1, a_2) \chi_{1,2}[ \phi (z_2;a'_2(1)) - \phi (z_2;a_2(1)) ]  \\
&= \lambda \frac {N-1}N  \int dz_1dz_2 f_2^N (z_1, z_2; I, S) \chi_{1,2}[ \phi (z_2;I) - \phi (z_2;S) ]\;. 
\end{split}
\end{equation*}
We see that we recover the same kinetic system obtained for the first model, provided the propagation of chaos  holds true. Actually this is not the case, at least for a suitable choice of parameters. Indeed 
macroscopic correlations in areas $O(R_0^2)$ could be created when the infected crowds do not have enough time to mix. We will discuss in Section \ref{sec:5} some numerical simulations supporting this observation.

We stress that having the same kinetic limit $N\to \infty$ for the two models introduced is not surprising, due to the separation of scales. The infection mechanism in the second model is much stronger as it involves a macroscopic portion of the population (instead of a single pair),
but the jumps in Model 1 have intensity $O(N)$ while in Model 2 they have intensity $O(1)$.

\subsection { SIR}

It is natural to compare the kinetic equations \eqref{B} with the well known SIR model for the evolution of the fraction of species $ \SS, \II, \RRR $ as a function of time: 
\be \label{SIR}
\begin{cases}
  \dot \SS = -\b \II \SS   \\
    \dot \II = \b \II \SS -\g \II \\
     \dot \RRR = \g \II \\
\end{cases}
\ee
for given $\g >0$ and $\b>0$. 

Define
$$
g(S;t)=M\, \SS(t), \,\, g(I;t)=M\, \II (t), \,\, g(R;t)=M\, \RRR(t),
$$
where $M = \frac 1 {2\pi |\L|}$.  Then  $g(A; t), A\in L$ solve \eqref{B} (as constant functions of $z$) provided that $\b= \frac {\lambda \pi R_0^2}{  |\L|}$.  
Therefore at equilibrium, namely when the distribution of each species is constant, the kinetic equations do not say more than the SIR model. However integrating Eq.s\,\eqref{B} with respect to $z$, we do not find closed equations for the fractions
$$A(t) :=\int dz f(z;A;t)\;, $$ which means that in case of inhomogeneous data the kinetic equations do provide a more detailed description. 

A more accurate SIR model takes into account also the possibility that the recovered agents may become susceptible after some time. The equations are:
\be \label{SIR2}
\begin{cases}
  \dot \SS = -\b \II \SS +\mu \RRR  \\
    \dot \II = \b \II \SS -\g \II \\
     \dot \RRR = \g \II -\mu \RRR \\
\end{cases}\;.
\ee
for $\mu >0$. The corresponding kinetic equations \eqref{B} modify slightly by joining the term $\pm \mu f(z;R)$ in the first and third equations respectively.

\section{Long time behaviour} \label{sec:3}
\setcounter{equation}{0}    
\def\theequation{3.\arabic{equation}}

Eq.s\,\eqref{SIR} have many stationary solutions, but a single one $ \left(\SS_\infty , \II_\infty ,  \RRR_\infty \right) $, which is the limit for $t\to \infty$ of the solutions to \eqref{SIR}. 
By the third equation of \eqref{SIR}
\be
\label{I}
\int_0^\infty \II (t) = \d <+\infty.
\ee
Setting now $\s=\SS+\II$ we have
$$
\dot \s=- \g \II  \,\,\,\, \text{and} \,\,\,\, \s(t)= \s_0 - \g \int_0^t d\t \II (\t)
$$
so that
$$
\s_\infty:= \lim_{t \to \infty} \s (t)= \s_0-\g \d.
$$
Setting $\AA_{\infty}:= \lim_ {t \to \infty} \AA (t)$, since $\SS (t)$ is decreasing, both $\II (t)$ and $\SS(t)$ are converging as $t\to \infty$ and, by \eqref{I}, $\II_\infty=0$. In conclusion

\be \label{asymt}
\begin{cases}
\SS_{\infty} =\s_0 -\g \d =\II_0+ \SS_0 -\g \d \\
\II_{\infty}  = 0 \\
\RRR_{\infty} =\RRR_0+ \g \d \\
\end{cases}\;.
\ee
Note that the dependence of the stationary solution on $\b$ is hidden in $\d$. 

A more precise determination of the asymptotic values is provided by the following well known argument. From
$$
\frac {d \SS}{d \RRR}= - \frac \b \g \SS,
$$
by using $\RRR_\infty+\SS_\infty=1$ and the assumption $\RRR_0=0$
one finds
$$
e^{ -\frac \b \g \SS_\infty} \frac \b \g \SS_\infty = \SS_0 \frac \b \g e^{- \frac \b \g }.
$$
Since $\max y e^{-y}=\frac 1 e$, given $\b$ and $\g$ one finds nonvanishing solutions for $\SS_\infty$.

\subsection{Time asymptotics of \eqref{B}}

We  try now a similar  analysis for the kinetic model.

We abbreviate by
$$
\NN_f (z; t)= f(z;S)(\chi_{R_0} * f(\cdot ;I))( z;t)
$$
the nonlinear terms in the r.h.s.\,of \eqref{B}, where $\chi_{R_0}=\chi ( |x| < R_0)$ is the indicator of the set $\{|x| < R_0\}$. We further set ${\bf L}:= -v\cdot \nabla_x +{\cal L}_1$, the generator of the random flight semigroup
$
U(t)=e^{{\bf L} t}\;.
$
Then we rewrite Eq.\,\eqref{B} 
\begin{equation}
\label{B1'}
\begin{cases}
\partial_t  f(S;t ) =  {\bf L}f(S;t)  - \lambda  \NN_f (t) \\
\partial_t  f(I;t)  ={\bf L}f(I;t)  +\lambda \NN_f (t) -\g f(I;t ) \\
\partial_t f(R;t) = {\bf L}f(R;t) + \g f (I;t) \\
\end{cases}
\end{equation}
and denote the initial data by $f_0(A)$.

We shall use that, for $h$ a probability density and arbitrary $t_0>0$,
\be
\label{ipoc}
\| U(t-t_0) h-M  \|_{L^1}  \leq 2 e^{-a(t-t_0)}
\ee 
(remind that $M=\frac 1 {2\pi |\L|}$) for some $a>0$.

Estimate \eqref{ipoc} is  well known in the framework of the very extended literature concerning the linear Boltzmann equation. Here, due to the simplicity of our equation,  we prefer to present a simple direct proof  in Appendix for the reader's convenience.

From the the third equation of \eqref{B1'} we have that
$$
f(R;t) =U(t) f_0(R) +\g \int_0^t d\t U(t-\t) f(I; \t)\;.
$$
Integrating both sides with respect to $z$, recalling that
$
A(t) =\int dz f(z;A;t)
$ and setting $A_0 = A(0)$,
 we obtain
$$
R(t) =  R_0 +\g \int_0^t d\t I(\t)\;.
$$
Clearly
$$
\int_0^\infty d\t I(\t) = \tilde \d < +\infty.
$$
Then setting
$$
f(S;t)+f(I;t) =\Sigma (t)
$$
we find
$$
\pa_t \Sigma (t)= {\bf L} \Sigma (t)-\g f(I;t) 
$$
so that 
$$
\frac d {dt} (S+I )=-\g I\;.
$$
Denoting by $A_{\infty} $ the asymptotic value of $A(t)$, proceeding as before we infer
(cf.\,\eqref{asymt})
\be 
\label{asymt1}
\begin{cases}
S_{\infty} = I_0+ S_0 -\g \tilde \d \\
I_{\infty}  = 0 \\
R_{\infty} =R_0+ \g \tilde \d \\
\end{cases}\;.
\ee
We cannot conclude that $\AA (t)=A(t)$ even if $\AA (0)=A(0)$ because, in general,  $\d \neq \tilde \d$.

Next, we study the asymptotic behaviour of the triple $\left(f(z;A;t)\right)_{A \in L}$ as $ t \to \infty$. Denote by $ g_\infty(A)=A_{\infty} M$ the equilibrium state. Setting 
\bea
&& \D(z;S)=f(z;S)-g_\infty(S)\;,\nn\\ 
&& \D(z;I)=f(z;I)-g_\infty(I)=f(z;I)\;,\nn\\
&& \D(z;R)=f(z;R)-g_\infty(R)\nn
\eea
we have, using that $g_\infty$ is a stationary solution to Eq.\,\eqref{B1'},
\begin{equation*}
\label{BD}
\begin{cases}
\partial_t  \D(S;t) ={\bf L} \D(S;t) - \lambda  [ \NN_f (t)-\NN_{g_\infty} (t)] \\
\partial_t \D(I;t)  = {\bf L} \D (I;t)  + \lambda [ \NN_f (t)-\NN_{g_\infty} (t)] -\g \D(I;t) \\
\partial_t \D (R;t) = {\bf L}\D ( R;t ) + \g \D (I;t) \\
\end{cases}\;;
\end{equation*}
hence
$$
\s (t) :=\D(S;t)+\D(I;t) =\D(S;t)+f(I;t)
$$
satisfies ($\s_0=\s(t_0)$)
\be
\s(t)= U(t-t_0) \s_0 + \g \int_{t_0}^t  d\t U(t-\t) \D(I;\t) \nn
\ee
for arbitrary $t_0>0$.

The two terms in this equation are controlled by
\be
\label{noto}
\g \left\| \int_{t_0}^t  d\t\, U(t-\t) \D(I; \t) \right\|_{L^1} \leq \g \int_{t_0}^t  d\t I (\t)\nn
\ee
and
\be
\label{ini}
\| U(t-t_0) \s_0 \|_{L^1}\leq  \| U(t-t_0) [ f(S;t_0)-S(t_0) M ]  \|_{L^1} + |S_{\infty} -S(t_0)| + I(t_0)\;.\nn
\ee
But, by \eqref{ipoc} and $\int \left(f(S;t_0)-S(t_0)\right) Mdz=0$,
 $$
  \| U(t-t_0) [ f(S;t_0)-S(t_0) M ]  \|_{L^1} 
  \leq 2 S(t_0) e^{-a(t-t_0)}  \leq  2  e^{-a(t-t_0)}\;.
 $$
Furthermore for any $\e >0$,
$$
\g \int_{t_0}^\infty  d\t I (\t) + |S_{\infty} -S(t_0)| + I(t_0) < \e
$$
provided that  $t_0$ is sufficiently large.
In conclusion
$$
\limsup_{t\to \infty} \| \s (t) \|_{L^1} < \e\;,
$$
thus
$$
\lim_{t \to \infty} \D(A;t) =0 \,\, \text{for} \,\, A \in L,
$$
in norm $L^1$.

Note that existence and uniqueness of the solutions to the system \eqref{B1'} follows by standard arguments, since the nonlinear term is Lipschitz continuous in $L^1$. Indeed by the identity
$$
\NN_f -\NN_h = \left(f-h \right)(S)\,\chi_{R_0}* f(I) +h(S)\,\chi_{R_0}* (f-h) (I) 
$$
and $ \| f (A;t) \|_{L^1} \leq 1$ it follows that
$$
\| \NN_f -\NN_h  \|_{L^1} \leq  C \| f-h \|_{L^1} \;, \qquad C>0\;.
$$

We summarize the discussion in the following statement.

\begin{prop}
Given an initial datum $ f_0(z;A) \geq 0, A \in L, z \in \G$, $f_0 \in L^1_{(x,v)}$ such that 
\be
\label{norm}
\sum_{A\in L} \int dz f_0( z;A)=1,
\ee 
  there exists a unique solution $f(z;A;t ) \geq 0$  to the system \eqref{B1'} for  any $t>0$, strongly differentiable in $L^1$ and preserving the normalization condition \eqref{norm}.
Moreover
\be
\lim_{t \to \infty} f(z;A;t) = g_\infty(A) \,\,\, \text {in} \,\,\, L^1
\ee  
where $g_\infty(A)=M A_\infty$, $M$ is the uniform normalized distribution on $\G$ and $A_\infty$ solves \eqref {asymt1}.
 \end{prop}

 \section{Particle approximation} \label{sec:4}
\setcounter{equation}{0}    
\def\theequation{4.\arabic{equation}}

In this section we derive the kinetic equations \eqref{B1'} starting from Model 1. This is a classical mean-field problem, which has been largely investigated in previous literature.  Observe that Model 1 is particularly simple because, once integrated over labels, the probability measure factorizes
$$
W^N_t(Z_N;t) := \sum_{A_N} W^N_t(Z_N;A_N) = (f(t) )^{\otimes N} (Z_N)
$$
and $f$ satisfies Eq.\,\eqref{rf}. We recall that we are assuming  \eqref{eq:fullindtz} and \eqref{eq:norm0} at time zero.
On the other hand the dynamics of labels generates correlations as
$$
f^N_2 ( z_1,z_2; a_1,a_2; t) \neq f^N_1 ( z_1;a_1; t) f^N_1 (z_2;a_2; t)\;.
$$
In order to show that such correlations are negligible as $N \to \infty$,
a straightforward method consists in establishing a hierarchy of equations for the marginals $f^N_j (Z_j;A_j;t)$. 
This is a purely algebraic computation, leading to the following result.
\begin{thm}
{ \it For any $t >0$ we have that the marginals of Model 1 are chaotic, and
$$
\lim_{N\to \infty} f^N_k (t)= f^{\otimes k}(t)
$$
in $L^1_{Z_k}$ and for any choice of $A_k$, where $ f=f(z;a;t)$ solves the kinetic system \eqref{B1'}. }
\end{thm}

\noindent
{\em Proof.}
We start by computing the adjoint of the generator \eqref{gen1}
\be
\label{gen*}
 {\cal L}^* = -{\cal L}_0+{\cal L}_1+ {\cal L}^*_d +{\cal L}^{N*}_{int}\;.\nn
\ee
The decay operator has adjoint
\be
{\cal L}^*_d W^N (Z_N; A_N) = \g \sum_{i=1} ^N 
 W^N ( Z_N; a_1, \ldots, \tilde a^*_i, \ldots, a_N) ( \d_{a_i,R}-\d_{a_i,I} ) \nn
\ee
where
$$
\tilde a^*_i=I \,\,\,\,\,\, \text {if} \,\,\,\,\,a_i=R; \qquad \tilde a^*_i=a_i  \,\,\,\,\, \text{otherwise}\;.
$$
The interacting generator is computed as follows.
Denoting $$ A^{i,j}_N (b,d)=\{ a_1, \dots, a_{i-1}, b, a_{i+1}, \dots, a_{j-1}, d, a_{j+1} ,\dots, a_N \} \qquad (i \neq j)$$
one has that
\bea
\label{dual}
&&\sum_{A_N} \int dZ_N  W^N (Z_N;A_N)  {\cal L}^N_{int} \Phi (Z_N;A_N)\nn \\
&& = \frac \lambda {N} \sum_{ i < j}\,\, \sum_{A_N} \int dZ_N W^N(Z_N;A_N) [ \Phi (Z_N ,A^{i,j}_N(a_i',a_j')) - \Phi (Z_N;A_N) ] \nn\\
&& = \frac \lambda {N} \sum_{ i < j}\,\, \sum_{A_N} \int dZ_N W^N(Z_N;A_N) \chi_{i,j} [ \d_{a_i, I}\d_{a_j,S} + \d_{a_i, S}\d_{a_j,I} ]
 \nn \\
&&\qquad\cdot  [ \sum_{ a_i', a_j'} \d_{a'_i, I}\d_{a'_j,I} \Phi (Z_N ,A^{i,j}_N(a_i',a_j') )-\Phi (Z_N ,A_N) ]\;, \nn 
\eea
from which we obtain
\bea
\label{dualgen}
{\cal L}^{N*}_{int}\, W^N(Z_N;A_N)&=& \frac \lambda {N} \sum_{ i < j}\,\, \chi_{i,j} \{ (W^N(Z_N;A^{i,j}_N (I,S))\\&&
+W^N(Z_N;A^{i,j}_N (S,I)))  \d_{a_{i, I }}  \d_{a_j, I } \nn \\
&&- W^N (Z_N;A_N) ( \d_{a_i, I}\d_{a_j,S} + \d_{a_i, S}\d_{a_j,I} )\}\;.\nn
\eea

The hierarchical equation for marginals is obtained (as for the well known BBGKY hierarchy) by computing the quantity
\be
\label{C}
\sum_{A'_{N-k}} \int dZ'_{N-k}{\cal L}^{N*}_{int}\, W^N (Z_k, Z'_{N-k} ;A_k, A'_{N-k})\;. \nn
\ee
We split the sum $ \sum_{ i < j}$ into three contributions. The first one for $i<j\leq k$ yields
$$
 \frac k N {\cal L}^{k*}_{int}\,  f^N_k\;.
$$
The second one for $k <i<j$ is vanishing. The third one, for $i\leq k, j>k$, gives, using the symmetry of $W^N$, 
\bea
&& \lambda \frac {N-k}N  \sum_{ i=1}^k \,\, \int dz_{k+1} \, \chi_{i,k+1} \nn\\
&&\qquad\cdot \Big\{ \d_{a_i,I}  \big( f_{k+1} ^N(Z_k,z_{k+1} ;A^{i,k+1}_{k+1} (I,S))
+f_{k+1} ^N(Z_k,z_{k+1} ;A^{i,k+1}_{k+1} (S,I)) \big)- \nn \\
&&\qquad\quad \big( \d_{a_i, I } f^N_{k+1}  (Z_k,z_{k+1} ;A_k, S) +  \d_{a_i,S} f^N_{k+1}  (Z_k,z_{k+1} ;A_k, I) \big)\Big\}  \nn \\
&&=   \lambda \frac {N-k}N  \sum_{ i=1}^k \,\, \int dz_{k+1} \chi_{i,k+1}
\\&&\qquad \quad \cdot \,\d_{a_i,S}  [  f^N_{k+1}  (Z_k,z_{k+1} ;A_{k+1}^{i,k+1}(I,I))
- f^N_{k+1}  (Z_k,z_{k+1} ;A_k,I)] \nn \\
&& =:
\lambda \frac {N-k}N C_{k+1} f^N_{k+1} (Z_k;A_k)\;.\nn
\eea
Here the last identity defines the hierarchical collision operator $C_{k+1}$.

In conclusion we find 
\be
\label{hie}
\pa_t f^N_k={\bf L} f^N_k + {\cal L}^*_d  f^N_k+ \frac k N {\cal L}^{k*}_{int}\,  f^N_k+ \lambda \frac {N-k}N C_{k+1} f^N_{k+1} 
\ee
for $k < N$. The last equation, for $k=N$, is nothing else than the  equation for the measure $W^N \equiv f^N_N$, that is
$$
\pa_t W^N= {\cal L}^* W^N\;.
$$

Coming back to the kinetic system \eqref{B1'}, we write it in the more compact form of a single equation
\be
\label{kineq}
\pa_t f(z;a;t)={\bf L}f(z;a;t) + {\cal L}^*_d\, f (z;a;t) + \lambda \,Q(f,f)(z;a;t)\nn
\ee
where
$$
Q(f,f)(z;a):= \int dz_1\, \chi (|x-x_1| <R_0)\, [ f(z;S) f(z_1;I) \d_{a,I} - f(z;a) f(z_1;I) \d_{a,S} ]\;.
$$
Consider now the sequence of products
$$
f_j(Z_j;A_j;t) := f^{\otimes j} (Z_j;A_j;t)\;.
$$
By direct inspection, we obtain
\be
\label{hieB}
\pa_t f_k={\bf L}f_k + {\cal L}^*_d \, f_k+ \lambda \,C_{k+1} f_{k+1} 
\ee
for $k \geq 1$.

We are now in position to conclude the proof following, for instance, the same strategy as for the inhomogeneous Kac model
(see e.g.\,\cite{PWZ}), which is inspired to the seminal paper by Lanford on the validity of the Boltzmann equation for hard sphere systems
\cite{Lanford}.  We remind the basic steps.

\begin{enumerate}
\item  The operator $C_{k+1}$ is controlled by
$$
\| C_{k+1} f_{k+1} \|_{L_{A}^1} \leq C k \| f_{k+1} \|_{L_A^1}
$$
for some $C>0$, 
where the norm $ \| \cdot \|_{L^1_A}$  is defined as
$$
\| f_k \|_{L^1_A} :=\sum_{A_k} \| f_k (A_k) \|_{L^1_{Z_k} }\;.
$$
\item We can represent the solutions of both hierarchies \eqref{hie} and \eqref{hieB} in terms of two series expansions which are converging in $L^1_A$ for short times, uniformly in $N$.
\item The Markov semigroup $ U^N_k(t) := e^{ \left({\bf L} +  \frac k N{\cal L}^{k*}_{int} \right)t}$ converges in $L^1_A$ in the limit $N\to \infty$ to 
$ U(t) = e^{ {\bf L} t}$, indeed 
$$
 \frac k N {\cal L}^{k*}_{int} =O\left( \frac {k^2}N\right)\;.
$$
Hence we have a term by term convergence of \eqref{hie} to \eqref{hieB}.
\item This allows us to achieve a short time convergence. But we have the a-priori estimate $ \| f_k^N(t) \|_{L^1_A} =1$, which allows us to iterate the procedure and reach arbitrary times. 
\end{enumerate}
\qed

 \section{Numerical simulations} \label{sec:5}
\setcounter{equation}{0}    
\def\theequation{5.\arabic{equation}}

We make use here of Monte Carlo method to  simulate the behaviour of Model $1$ and Model $2$, and compare the evolution of the population  fractions $S(t)$, $I(t)$ and $R(t)$ 
 with the solution of the SIR equations \eqref{SIR}.

Let us describe the setting of the particle simulation.
The spatial domain is the torus $\L=(0,D)^2$ with $D=500$.  
At time $0$ we consider $N$ particles uniformly distributed in space with  uniformly distributed velocities in ${\mathbb S}$, so that the gas as a whole is at equilibrium.
We focus on two different initial distributions of labels. In the first  case, a fraction $I(0)$ of infected agents are labeled as $I$, and these particles are chosen uniformly. In the second case, the $I$ particles at time $0$ are all the particles lying in a disk of area $I(0) |\Lambda|$. All the remaining agents are susceptible, hence we fix $R(0)=0$.
In the following we will refer to these two initial distributions as homogeneous and concentrated initial data, respectively. In all the experiments reported below, we set $I(0)=\frac \pi {100}$, $S(0)=1-I(0)$, $R(0)=0$.

We recall here that the kinetic equations \eqref{B} reduce to the SIR--model \eqref{SIR} for uniform data,  with  the corresponding parameter $\beta$ to be chosen as $\beta=\frac{\lambda \pi R_0^2}{|\Lambda|}$. Also, the SIR--model asymptotics, once fixed the initial data, depends only on the ratio $\frac{\beta}{\gamma}$. For fixed ratio $\frac{\beta}{\gamma}$, the actual values of $\beta$ and $\gamma$ influence only the time scale of the evolution. 
As we shall see, inhomogeneous initial data can instead modify considerably the evolution of the population fractions.

\begin{figure}[th!]
\begin{picture}(200,145)(0,0)
\put(-5,-3){
  \includegraphics[width=0.55\textwidth]{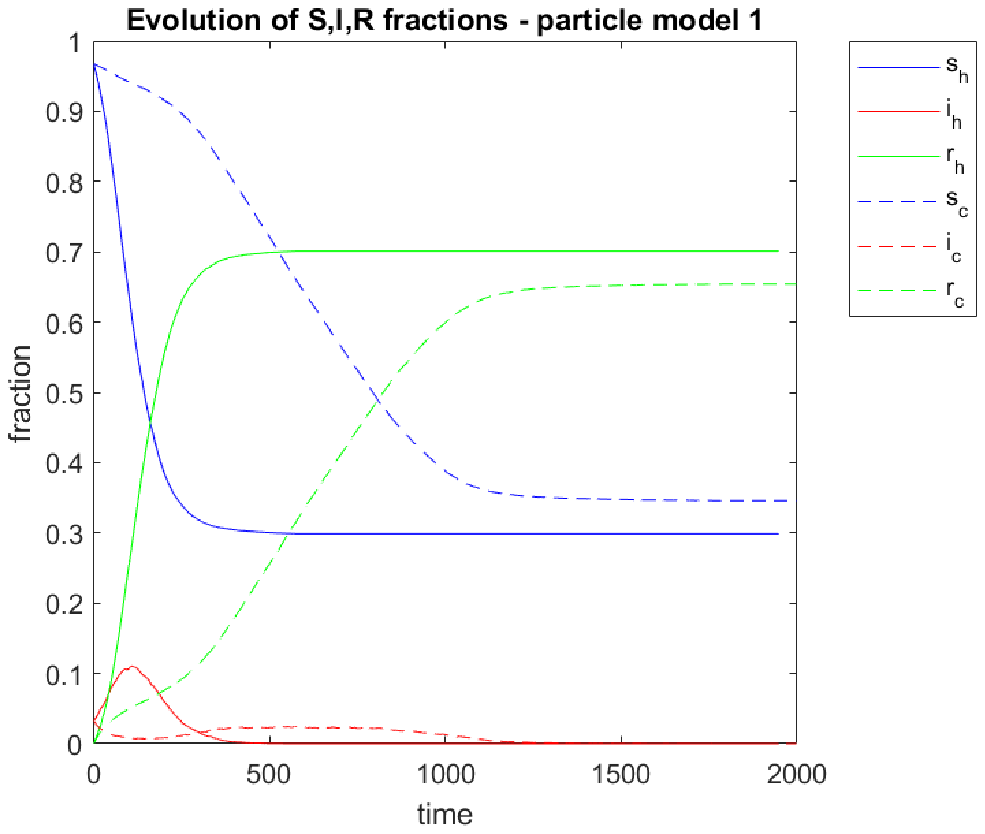}
}
\put(194,0){
  \includegraphics[width=0.385\textwidth]{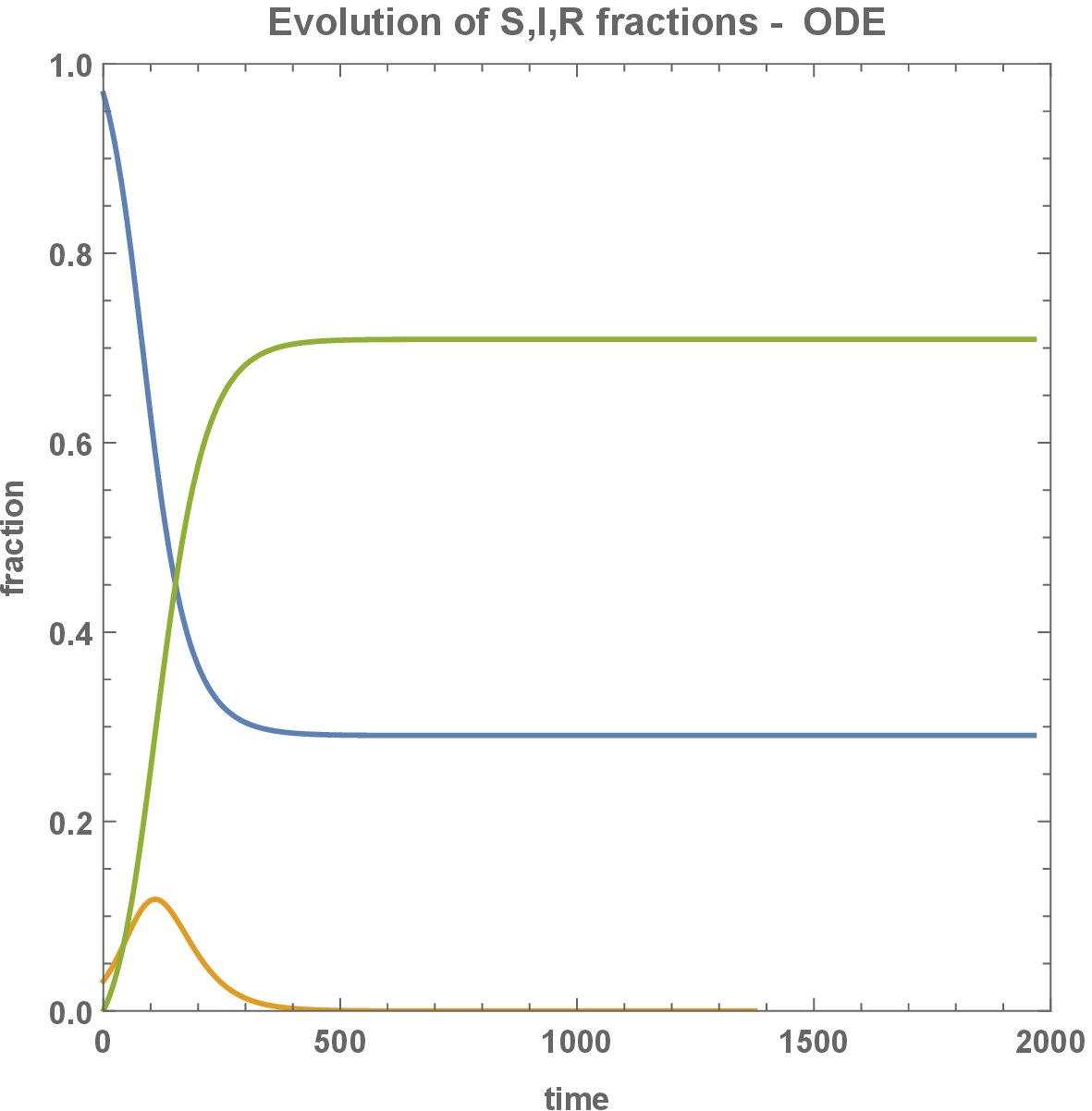}
}
\put(331,104){
  \includegraphics[width=0.05\textwidth]{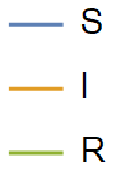}
  }
\end{picture}
\caption{ 
Parameters:  $\lambda=20$, $R_0=15$, $\gamma=1/30$.
\\Left panel: evolution of $S$, $I$, $R$  fractions for a particle system evolving according to Model $1$. 
The solid and dashed lines represent, respectively, the case of homogeneous and concentrated initial data. $N=50000$.
\\Right panel: numerical solution of \eqref{SIR} 
with $\beta=\frac{\lambda \pi R_0^2}{|\Lambda|}$.
}
\label{f:001}
\end{figure}

We consider first the dynamics  of Model $1$,
with parameters $\lambda=20$, $R_0=15$, $\gamma=1/30$; see Fig.\,\ref{f:001}. 
The result verifies the correspondence between homogeneous particle model and SIR--model, and the different behaviour in the case of concentrated initial data.
In  the latter case, the spread of infected particles is much slower, implying that the infected population reaches a lower peak and in a longer time. 
The asymptotic values $S_\infty$ and $R_\infty$ are also affected. The quantity $\tilde{\delta}$ appearing in
\eqref{asymt1} is indeed different in the homogeneous 
and in the concentrated case.

When considering Model $2$ with the same parameters as in Figure $1$, the evolution is very far from the SIR--model, even for homogeneous initial data, and only a small fraction of susceptible agents is infected before the extinction of the infected population.
Note that the evolution of Model $2$, as pointed out in Section \ref{ssect:mod2}, only gives rise to the kinetic description when the propagation of chaos holds true in the limit $N\to\infty$, which might not hold for this choice of parameters.

Let us now fix the ratio $\frac{\lambda \pi R_0^2}{|\Lambda|} / \g$ and the domain $|\Lambda|$, and let us initially fix also the value of $\gamma$.
We can choose different values for $R_0$ and $\lambda$  and  try to identify the regime for which the homogeneous particle system behaves as 
the SIR-ODE model. We find that by decreasing $R_0$ and increasing  $\lambda$ accordingly, keeping  the product $\lambda  R_0^2$ fixed, the particle system approaches the SIR-ODE model behaviour, while it is far from it for large $R_0$. 
 The rate of approach also depends on $\gamma$ (for smaller $\gamma$, higher value of $R_0$ is required). However, this is true up to a certain threshold. Indeed, when the value for $R_0$ is too large, no choice of $\gamma$ and $\lambda$ can work (think of $R_0>D$).

The spread of infected particles is favoured when agents of type $I$ are surrounded by a large number of susceptible. 
This does not happen in general when the labels have strongly non--homogeneous distribution.
The dynamics in Model $2$ produces such inhomogeneities  in disks of radius $R_0$.
However, if the decay rate $\gamma$ is sufficiently small, infected particles have sufficient time to mix with other agents in the surrounding space, before becoming recovered. 
Instead for $R_0$ large, the infected agents are unable to exit the shielded region.

\begin{figure}[htp!]
\begin{picture}(200,145)(0,0)
\put(-6,-3){
  \includegraphics[width=0.55\textwidth]{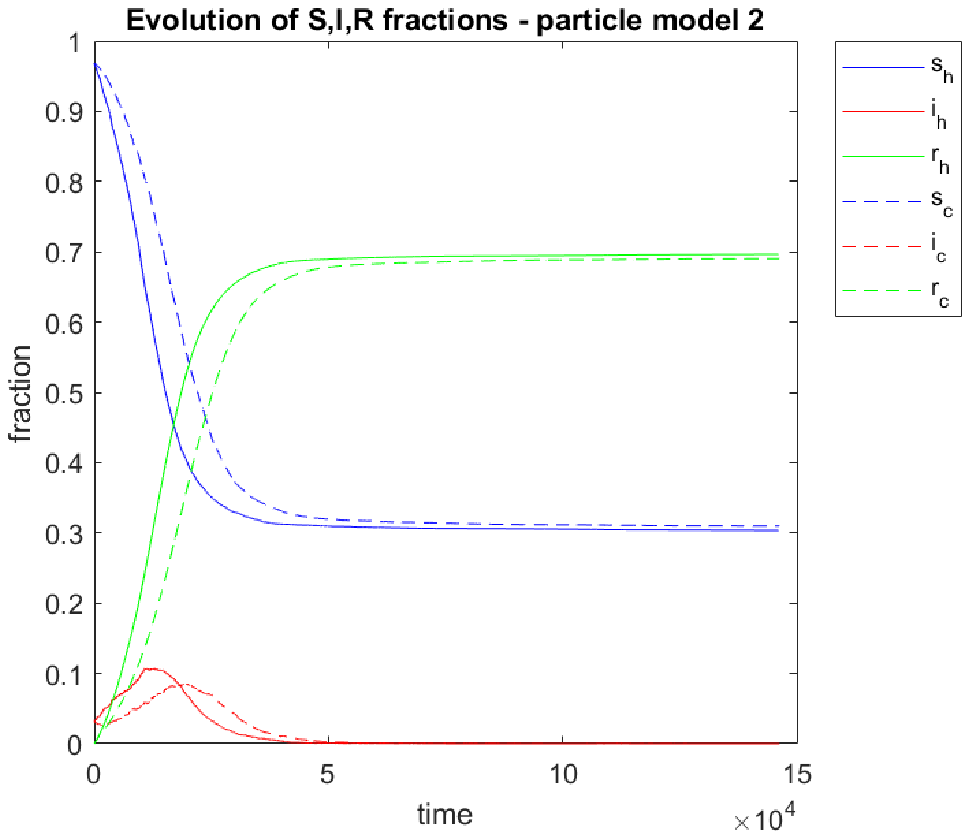}
}
\put(194,0){
  \includegraphics[width=0.374\textwidth]{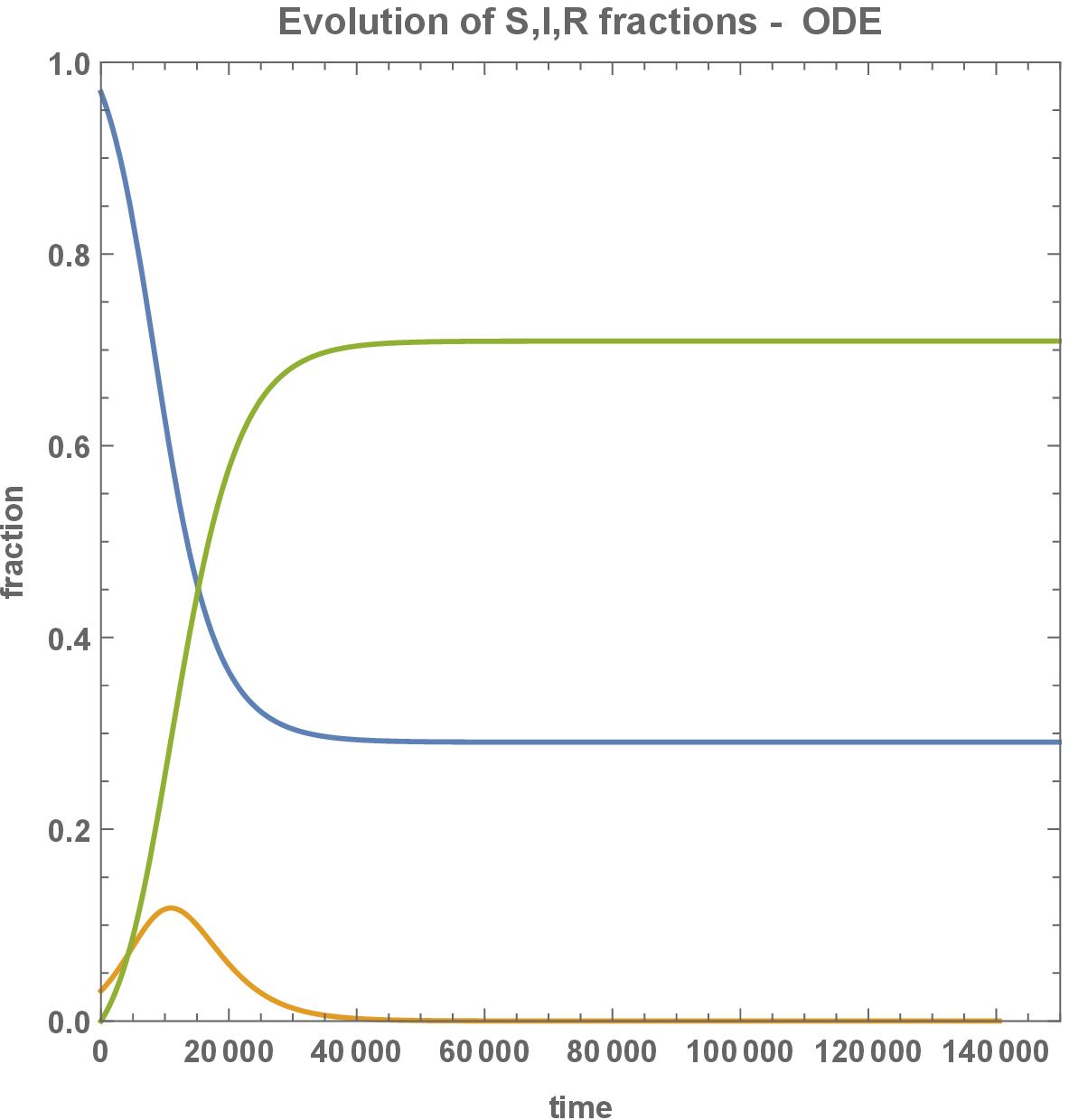}
}
\put(331,104){
  \includegraphics[width=0.05\textwidth]{legenda.eps}
  }
\end{picture}
\caption{ 
Parameters:  $\lambda=2$, $R_0={15}/\sqrt{10}$, $\gamma=1/3000$.
\\Left panel: evolution of $S$, $I$, $R$  fractions for a particle system evolving according to Model $2$. 
The solid and dashed lines represent, respectively, the case of homogeneous and concentrated initial population of $I$ agents, $N=200000$.
\\Right panel: numerical solution of \eqref{SIR} 
with $\beta=\frac{\lambda \pi R_0^2}{|\Lambda|}$.
}
\label{f:002}
\end{figure}

\begin{figure}[htp!]
\begin{picture}(200,145)(0,0)
\put(-6,-3){
  \includegraphics[width=0.55\textwidth]{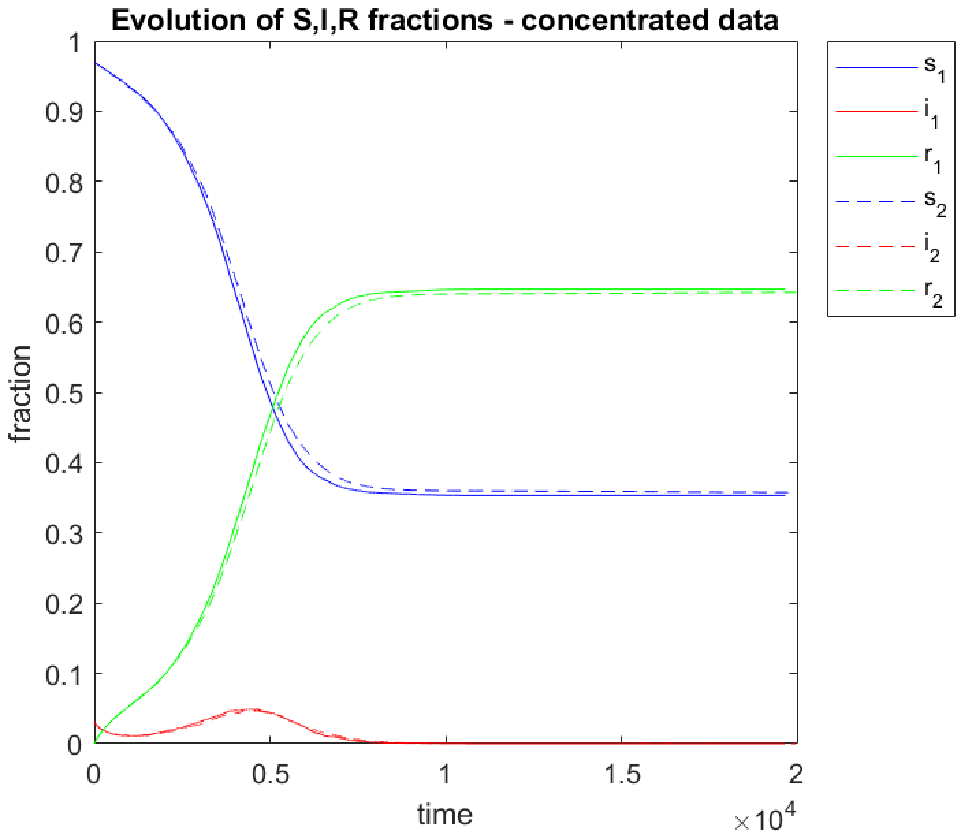}
}
\put(194,0){
  \includegraphics[width=0.374\textwidth]{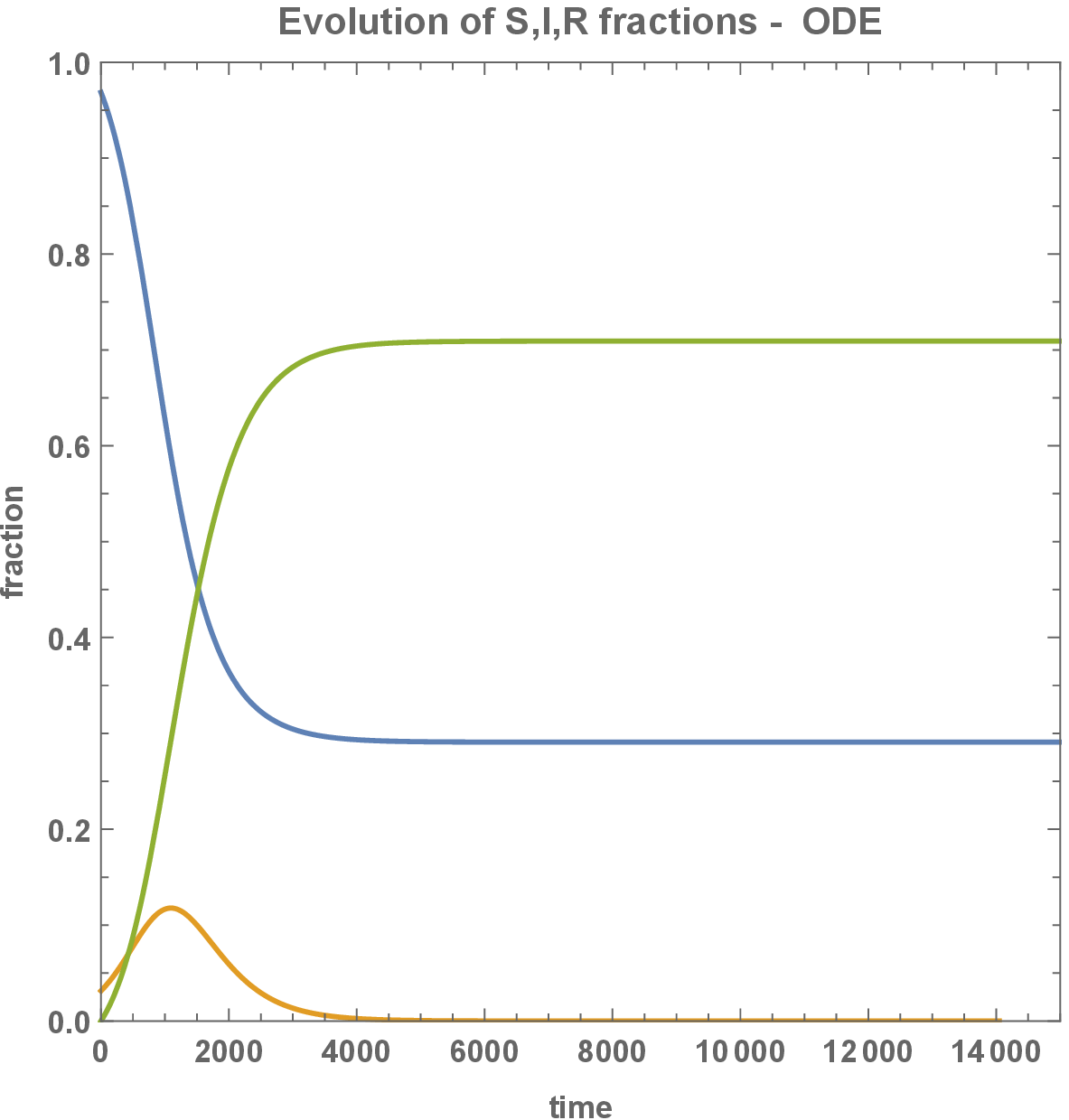}
}
\put(331,104){
  \includegraphics[width=0.05\textwidth]{legenda.eps}
  }
\end{picture}
\caption{ 
Parameters:  $\lambda=200$, $R_0=15/{10}$, $\gamma=1/300$.
\\Left panel: evolution of $S$, $I$, $R$  fractions for a particle system evolving according to Model $1$ (solid lines) and to Model $2$ (dashed lines), for concentrated initial population infected agents. $N=60000$ and $N=200000$ respectively.
\\Right panel: numerical solution of \eqref{SIR} 
with $\beta=\frac{\lambda \pi R_0^2}{|\Lambda|}$.
}
\label{f:003}
\end{figure}

We see in Fig. \ref{f:002} a second set of simulations where the parameters are such that Models $2$ and SIR Equations match (keep the same ratio $\frac{\beta}{\gamma}$).
Now, since $\gamma$ and $\lambda$ are small, particles have enough time to mix before being involved in a new infection.
Moreover for the same reason, the particles reach a homogeneous distribution quickly, so that the difference homogeneous and concentrated initial datum is not so significant in this case. 
In Fig.\,\ref{f:003}, we propose a last set of simulations, for an intermediate regime (and same $\frac{\beta}{\gamma}$).

\newpage
 \section{Appendix. Proof of estimate \eqref{ipoc}}\label{sec:A}
\setcounter{equation}{0}    
\def\theequation{A.\arabic{equation}}

We fix $t_0=0$ with no loss of generality.
Consider the following equation for a probability density $f(x,v,t)$ with $x\in \L$ and $v \in {\mathbb S}$:
\begin{equation*}
\label{lin}
\begin{cases}
\partial_t f +v\cdot \nabla_x f= \frac {\rho_f}{2\pi}  - f  \\
f(x,v,0)=f_0(x,v) \\
\end{cases}\; 
\end{equation*}
where $\rho_f (x)=\int_{{\mathbb S}} dw f(x,w) $ is the spatial density associated to $f$. We can write the solution explicitly, by iteration of the Duhamel formula:
\bea
f(x,v,t)&=& e^{-t} \sum_{n\geq 0} \left(\frac 1{2\pi}\right)^n \int_0^t dt_1 \cdots \int_0^{t_{n-1}} dt_n \int dw_1 \cdots \int dw_n \nn \\
&& f^p_0 (x-v(t-t_1)-w_1(t_1-t_2)\cdots -w_n t_n, w_n)\;,\nn
\eea
where $f^p_0$ defined in the whole $\R^2 \times {\mathbb S}$, is the periodic extension of $f_0$ from the square $\L$.

Let $P_t(z_0,z)$ be the transition probability from $z_0=(x_0,v_0)$ to $z=(x,v)$ in time $t>0$. 
 For $y \in \R^2$ and $z_0\in \L$, we introduce
$
\d^p(z-z_0):= \sum_{\pi} \d ( x- (x_0)_{\pi},v-v_0)
$
where $\d$ is the Dirac delta, $\pi =(k_1,k_2)$ is a pair of integers and $$(x_0)_{\pi}=((x_0)_{1}+k_1D, (x_0)_{2}+k_2 D)$$ are periodic images of $x_0$. We also denote abusively $\d^p(x-x_0):= \sum_{\pi} \d ( x- (x_0)_{\pi})$.
Then
\bea
\label{pr}
P_t(z_0,z) &=& e^{-t}\, \d^p (x-vt-x_0,v-v_0) + e^{-t} \sum_{n\geq 1} \left(\frac 1{2\pi} \right)^n \int_0^t dt_1 \cdots \int_0^{t_{n-1}} dt_n \nn\\ && \int dw_1 \cdots \int dw_{n-1} \,
 \d^p  (x-v(t-t_1)-w_1(t_1-t_2)\cdots -v_0 t_n -x_0)\;. \nn
\eea

$P_t$  is not absolutely continuous and we are looking for a lower bound, uniform in $z$ and $z_0$. The above formula is a series of positive terms which we call $P^{(n)}_t(z_0,z)$. They are absolutely continuous for $n \geq 2$. We focus then on the simplest contribution $n=2$
\be
\label{n=2}
P^{(2)}_t(z_0,z):= e^{-t}\!  \left(\frac 1{2\pi} \right)^2\!\! \int_0^t \!dt_1 \! \int_0^{t_1} dt_2 \int dw \,\,\, \d^p( x-v(t-t_1)-w(t_1-t_2) -v_0t_2- x_0)\;.\nn
\ee
Changing to the variable
$
\xi = w (t_1-t_2)
$
one finds that $P^{(2)}_t(z_0,z) $ is equal to
\begin{equation*}
\begin{split}
\label{n=2'}
&e^{-t}   \left(\frac 1{2\pi} \right)^2 \int_0^t dt_2\int d\xi \,\,  \chi( |\xi| \leq t-t_2) \frac 1 {|\xi|} \nn  \, 
\d^p( x-v(t-t_2) - \xi + |\xi | v - v_0t_2- x_0 ) \\
&\geq  e^{-t}  \left(\frac 1{2\pi} \right)^2\! \int_0^{\frac t2}  dt_2\!  \int d\xi \,\,  \chi( |\xi| \leq \frac t2)  
\, \d^p( x-v(t-t_2)- \xi + |\xi | v -v_0t_2- x_0 )  \,\, \frac 2 {t} \,.
\end{split}
\end{equation*}
If $t$ is large enough, $\xi$ spans at least a  square in the two-dimensional lattice of side $D$ and hence the above integral is not vanishing. Moreover the Jacobian $ | \frac {\pa \eta } {\pa \xi } |$  of the transformation
$$
\xi- |\xi | v  \to \eta
$$ 
is
$$
1 - \frac {v_2 \xi_2} {|\xi|} -  \frac {v_1 \xi_1} {|\xi|}
$$
with inverse bounded from below by
$$
\frac {|\xi|}{ |\xi | + |v_1| |\xi_1| + |v_2| |\xi_2|} \geq \frac 1 { 1+|v_1|+|v_2| } \geq \frac 13\;.
$$
Therefore the last integral in $d\xi$ is bounded from below by $1$, which implies
\be
\label{lb}
P_t^{ac}(z_0,z) \geq C e^{-t}
\ee
for some $C>0$ (independent  of  $z,z_0$), provided that $t$ is large enough. 
Here $P_t^{ac}$ is the absolutely continuous part of $P_t$.

From \eqref{lb} we shall prove that
\be
\label{fin}
\| P_t (z_0, \cdot ) - P_t (z_1, \cdot ) \|_{TV} \leq 2 \r
\ee
with $\r <1$ and this is enough (see e.g.\,\cite{Kul}) to conclude that 
\be
\label{fin1}
\| P_{nt} (z_0, \cdot ) - M \|_{TV} \leq 2 \r^n\,\nn
\ee
from which in turn we obtain for all $t>0$
$$
\| U(t) h-M  \|_{L^1}  \leq 2 e^{-at}\;,
$$ 
for some (certainly not optimal) $a>0$.

To prove \eqref{fin}, we introduce the Wasserstein distance with the discrete metric
$d(z,z')=1$ if $z \neq z'$,  $d(z,z)=0$:
$$
\WW ( \mu, \nu) = \inf _{R \in C(\mu,\nu)} \int_{\Gamma\times\Gamma} d(z,z')\,R(dz,dz') 
$$
where $C(\mu,\nu)$ is the family of  couplings between the probability measures $\mu$ and $\nu$ ($R$ is a measure on the product space having $\mu$ and $\nu$ as marginals).
We have that (see for instance \cite {Vill} Eq.\,(13), p.7)
$$
\| \mu -\nu \|_{TV} = 2  \,\, \WW (\mu;\nu)\;.
$$
To control $\WW (P_0,P_1)$ with $P_i(dz)=P_t(z_i, dz) $, $i=0,1$, we introduce the following explicit $R_0 \in C(P_0,P_1)$:
$$
R_0( dz,dz') = \d (z-z') \lambda (z)dz\,dz' + \frac { (P_0(dz)- \lambda (z)dz )(P_1(dz')- \lambda (z')dz' ) } {1-\int dz \lambda (z) }
$$
where $\lambda (z) = \min ( P_0^{ac} (z), P_1^{ac} (z))$ and $P_i^{ac}$ is the density of the absolutely continuous part of 
$P_i$.
By \eqref{lb} we obtain that
$$
\WW (P_1, P_2) \leq \int  d(z,z')\,R_0 (dz; dz') \leq 1-\int dz \lambda (z) \leq 1- C e^{-t} =:\rho.
$$
This concludes the proof. \qed

\subsection*{Acknowledgement}
We are indebted to Nadim Sah for pointing out the relevance of superspread processes, which motivated part of this work.

\end{document}